\documentclass[11pt,reqno]{amsart}
\usepackage{amsmath,amsthm,amssymb,amscd,url,enumerate}
\usepackage{graphicx}
\usepackage{afterpage}
\usepackage{tikz}
\usepackage[all]{xy}
\usepackage[colorlinks=true]{hyperref}
\usepackage[width=5.8in,height=8.5in, centering]{geometry}

\newcommand{\TITLE}{Visualising the arithmetic of imaginary quadratic fields}
\newcommand{\TITLERUNNING}{Visualising the arithmetic of imaginary quadratic fields}
\newcommand{\DATE}{\today}

\theoremstyle{plain} 
\newtheorem{theorem}{Theorem} 

\newtheorem{proposition}[theorem]{Proposition}
\newtheorem{lemma}[theorem]{Lemma}
\newtheorem{corollary}[theorem]{Corollary}

\theoremstyle{definition}
\newtheorem{definition}[theorem]{Definition}

\theoremstyle{remark}

\numberwithin{theorem}{section}

  {\end{list}}

  {\end{list}}

\newcommand{\tightoverset}[2]{%
  \mathop{#2}\limits^{\vbox to -.5ex{\kern-1.05ex\hbox{$#1$}\vss}}}

\def\Pcal{{\mathcal P}}

\def\Scal{{\mathcal S}}

\newcommand{\CC}{\mathbb{C}}

\newcommand{\HH}{\mathbb{H}}

\newcommand{\QQ}{\mathbb{Q}}
\newcommand{\RR}{\mathbb{R}}
\newcommand{\MM}{\mathbb{M}}
\newcommand{\ZZ}{\mathbb{Z}}

\newcommand{\MOD}[1]{~(\textup{mod}~#1)}
\renewcommand{\pmod}{\MOD}

\renewcommand{\setminus}{\smallsetminus}

\newcommand\PSL{\operatorname{PSL}}
\newcommand\SL{\operatorname{SL}}
\newcommand\PGL{\operatorname{PGL}}
\newcommand\OK{\mathcal{O}_K}
\newcommand\Of{\mathcal{O}_f}
\newcommand\Picf{\mathcal{P}{ic}(\Of)}

\newcommand{\OSK}{\widehat{\Scal}_K}
\newcommand{\SK}{{\Scal}_K}

\newcommand{\Mobplus}{\text{M\"ob}_+(2)}

\newcommand\sensual{AppPrevious}

\newcommand\ntone{MR1971245}

\newcommand\gttwo{MR2183489}
\newcommand\gtone{MR2173929}

\newcommand\Cox{MR1028322}

\newcommand\Neukirch{MR1697859}
\newcommand\CohnBook{MR594936}

\newcommand\BourgainKontorovich{MR3211042}
\newcommand\FuchsBulletin{MR3020827}

\newcommand\Cohn{MR0207856}
\newcommand\Unreasonable{MR2805033}
\newcommand\Kocik{kocikminkowski}

\newcommand\Schmidt{MR0422168}
\newcommand\SchmidtEisenstein{MR698164}
\newcommand\SchmidtEleven{MR0485715}
\newcommand\SchmidtTwo{MR2811562}
\newcommand\SchmidtFarey{MR0245525}

\newcommand\Bianchi{MR1510727}

\newcommand\Serre{MR0272790}
\newcommand\FineAmalgam{MR920977}
\newcommand\FineBook{MR1010229}
\newcommand\Lubotzky{MR673125}
\newcommand\otherpaper{otherpaper}
\newcommand\Northshield{MR2250043}
\newcommand\BourgainFuchs{MR2813334}

\title[\TITLERUNNING]{\vspace*{-1.3cm} \TITLE}

\author{Katherine E. Stange}

\date{\DATE}
\address{%
Department of Mathematics, University of Colorado,
Campux Box 395, Boulder, Colorado 80309-0395}
\email{kstange@math.colorado.edu}

\thanks{Project sponsored by the National Security Agency under Grant Number H98230-14-1-0106.  The United States Government is authorized to reproduce and distribute reprints notwithstanding any copyright notation herein.
}

\begin{document}


\begin{abstract}
        We study the orbit of $\widehat{\RR}$ under the Bianchi group $\PSL_2(\OK)$, where $K$ is an imaginary quadratic field.  The orbit, called a Schmidt arrangement $\SK$, is a geometric realisation, as an intricate circle packing, of the arithmetic of $K$.  This paper presents several examples of this phenomenon.  First, we show that the curvatures of the circles are integer multiples of $\sqrt{-\Delta}$ and describe the curvatures of tangent circles in terms of the norm form of $\OK$.  Second, we show that the circles themselves are in bijection with certain ideal classes in orders of $\OK$, the conductor being a certain multiple of the curvature.  This allows us to count circles with class numbers.  Third, we show that the arrangement of circles is connected if and only if $\OK$ is Euclidean.  These results are meant as foundational for a study of a new class of thin groups generalising Apollonian groups, in a companion paper.
\end{abstract}

\maketitle

\section{Introduction}

{\bf Schmidt arrangements.}
The matrix groups $\PSL_2(\OK)$ are called Bianchi groups when $\OK$ is the ring of integers of an imaginary quadratic field $K$.  They are named for Bianchi's initial study in the 1890's \cite{\Bianchi}, as a natural generalisation of the modular group, $\PSL_2(\ZZ)$, and as a natural family of discrete subgroups of $\PSL_2(\CC)$, which can be realised as the group of isometries of hyperbolic $3$-space.  

The Bianchi groups have a beautiful visual manifestation as M\"obius transformations of the extended complex plane, $\widehat{\CC} = \CC \cup \{ \infty\}$.  We will identify the group $\Mobplus$ of orientation-preserving M\"obius transformations acting on $\widehat{\CC}$ with $\PSL_2(\CC)$, via the association
\[
        \left( z \mapsto \frac{\alpha z + \gamma}{\beta z + \delta} \right)
\leftrightarrow
\begin{pmatrix} \alpha & \gamma \\ \beta & \delta \end{pmatrix}.
\]
M\"obius transformations act transitively on the collection of circles in $\widehat{\CC}$ (where straight lines are circles through $\infty$).  

\begin{figure}
        \includegraphics[height=3in]{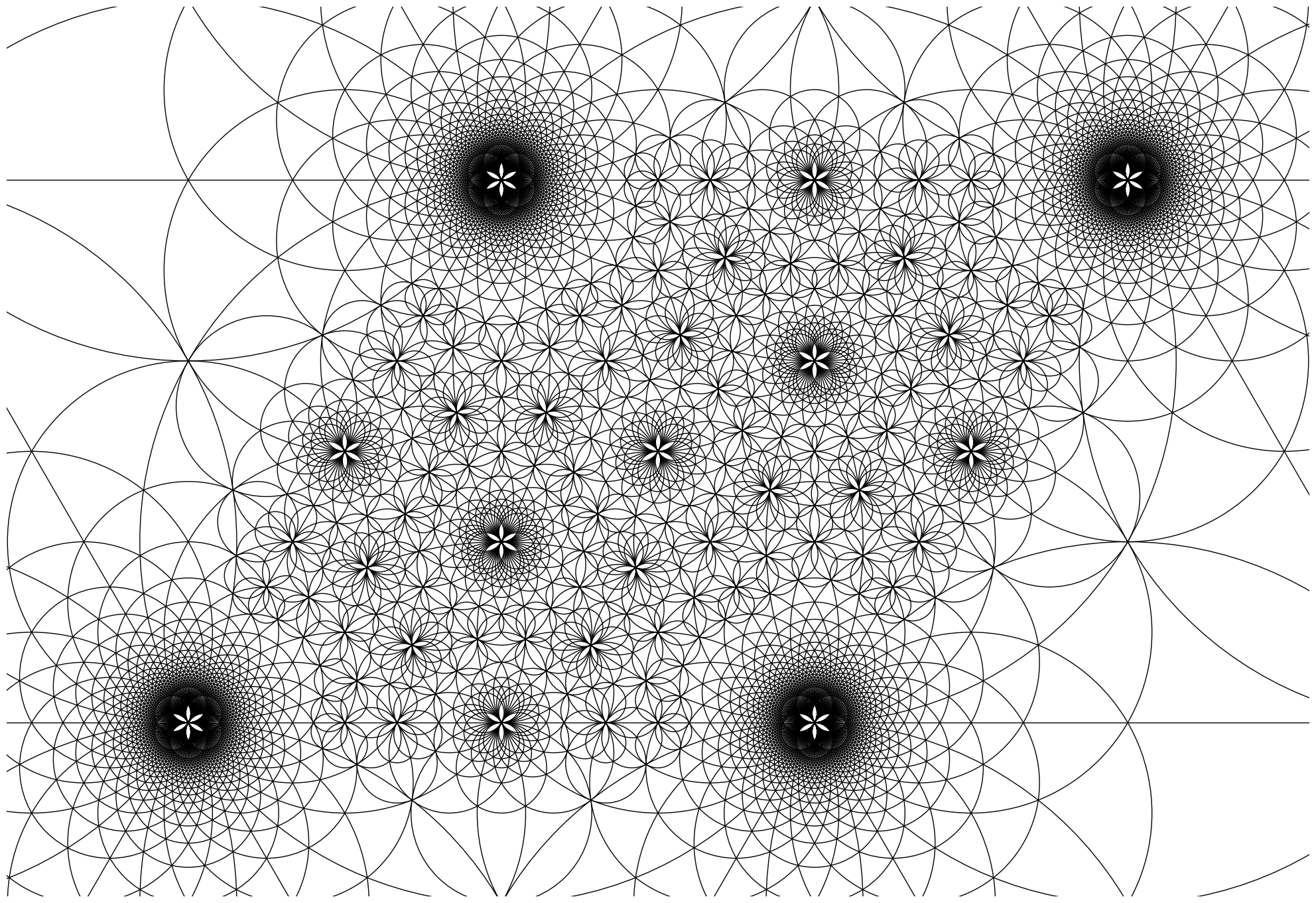} 
        \caption{Schmidt arrangement $\Scal_{\QQ(\sqrt{-3})}$.
        The image includes those circles of curvature bounded by $20$ intersecting the fundamental parallelogram of the ring of integers or its boundary.}
\label{fig:exampleSK3}
\end{figure}

\begin{figure}
        \includegraphics[height=3in]{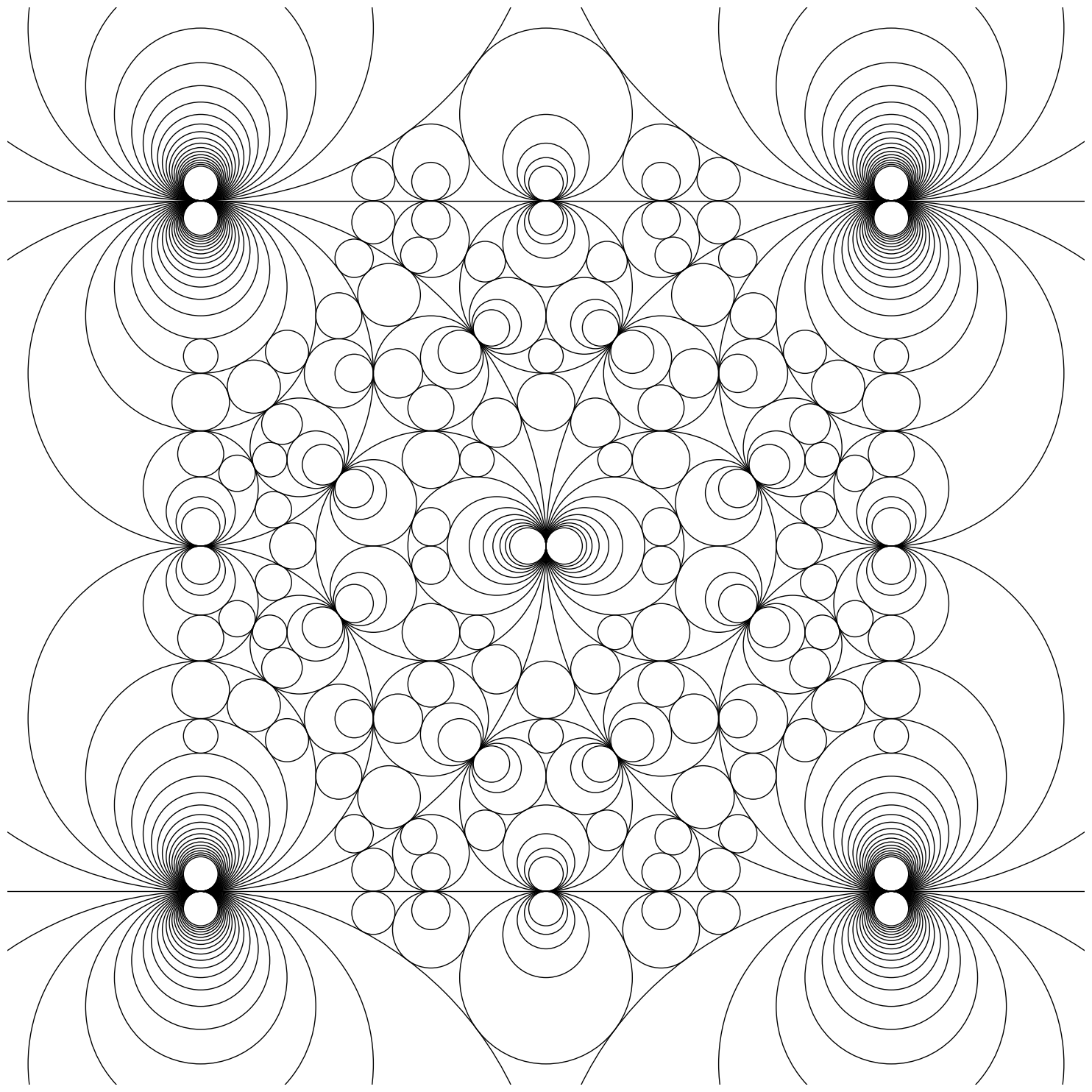}
        \caption{Schmidt arrangements ${\Scal}_{\QQ(\sqrt{-1})}$.  
        The image includes those circles of curvature bounded by $20$ intersecting the fundamental parallelogram of the ring of integers or its boundary.}
\label{fig:exampleSK}
\end{figure}

\begin{figure}
        \includegraphics[height=2.7in]{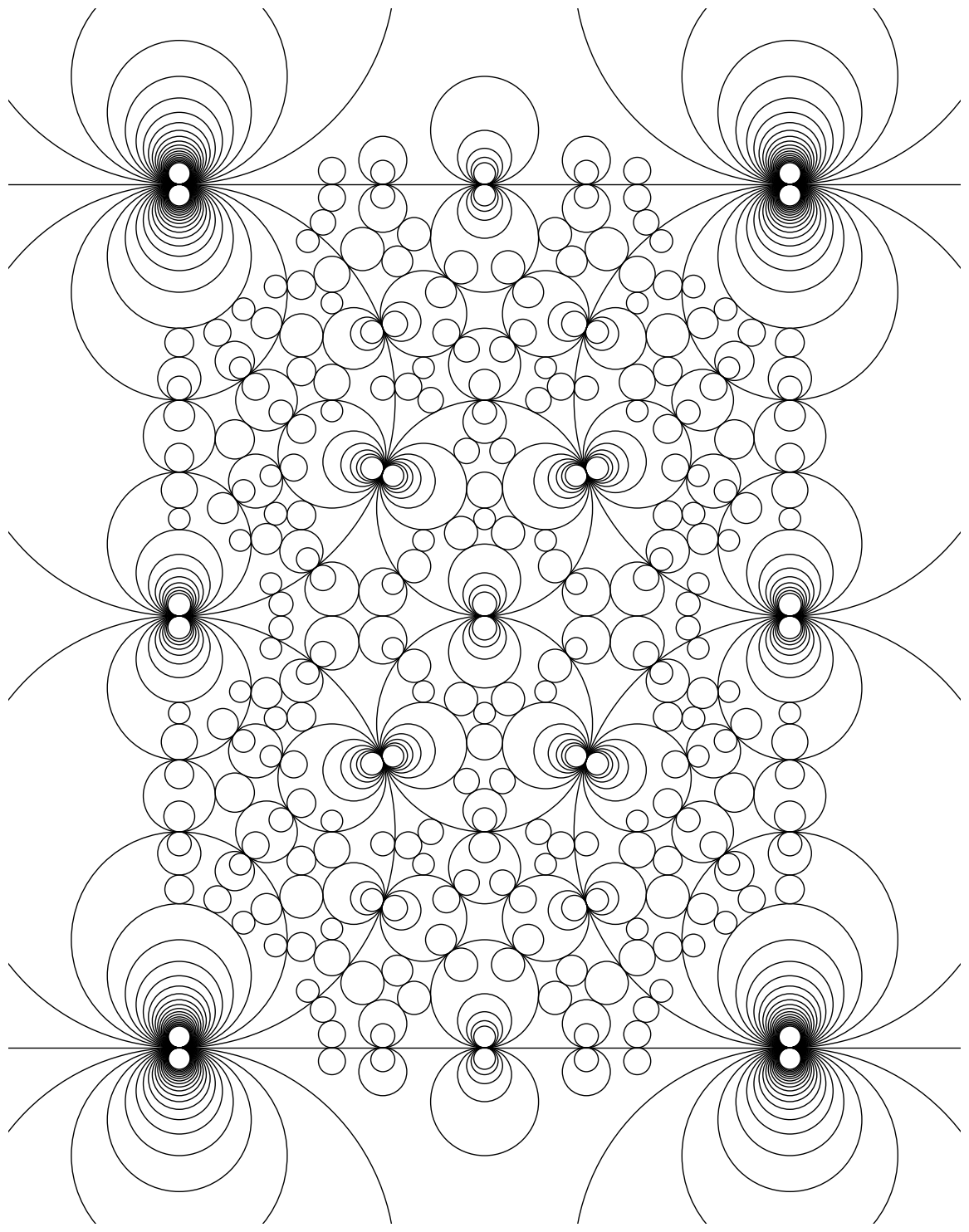} \quad
        \includegraphics[height=2.7in]{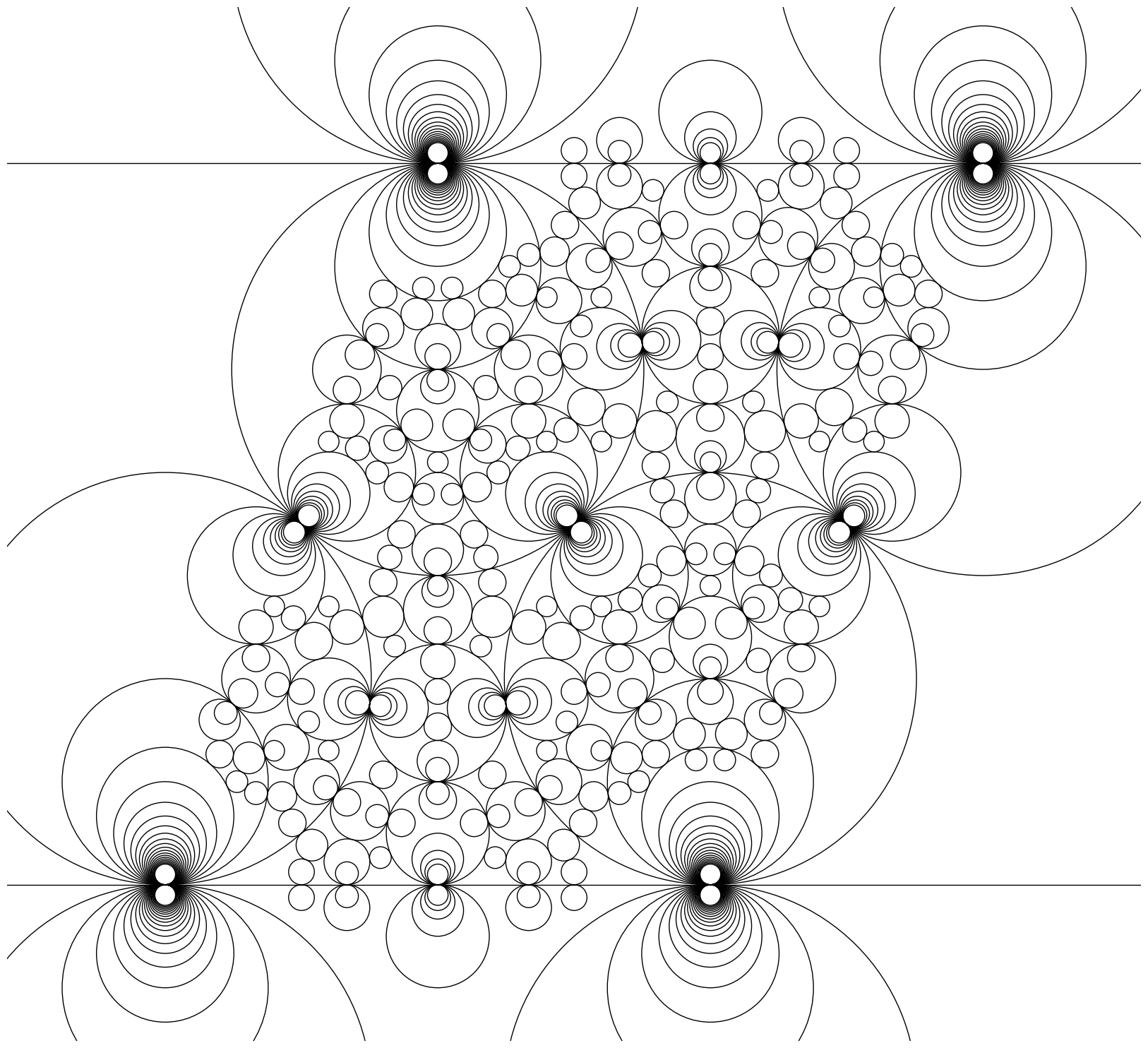} \\
        \includegraphics[height=2.9in]{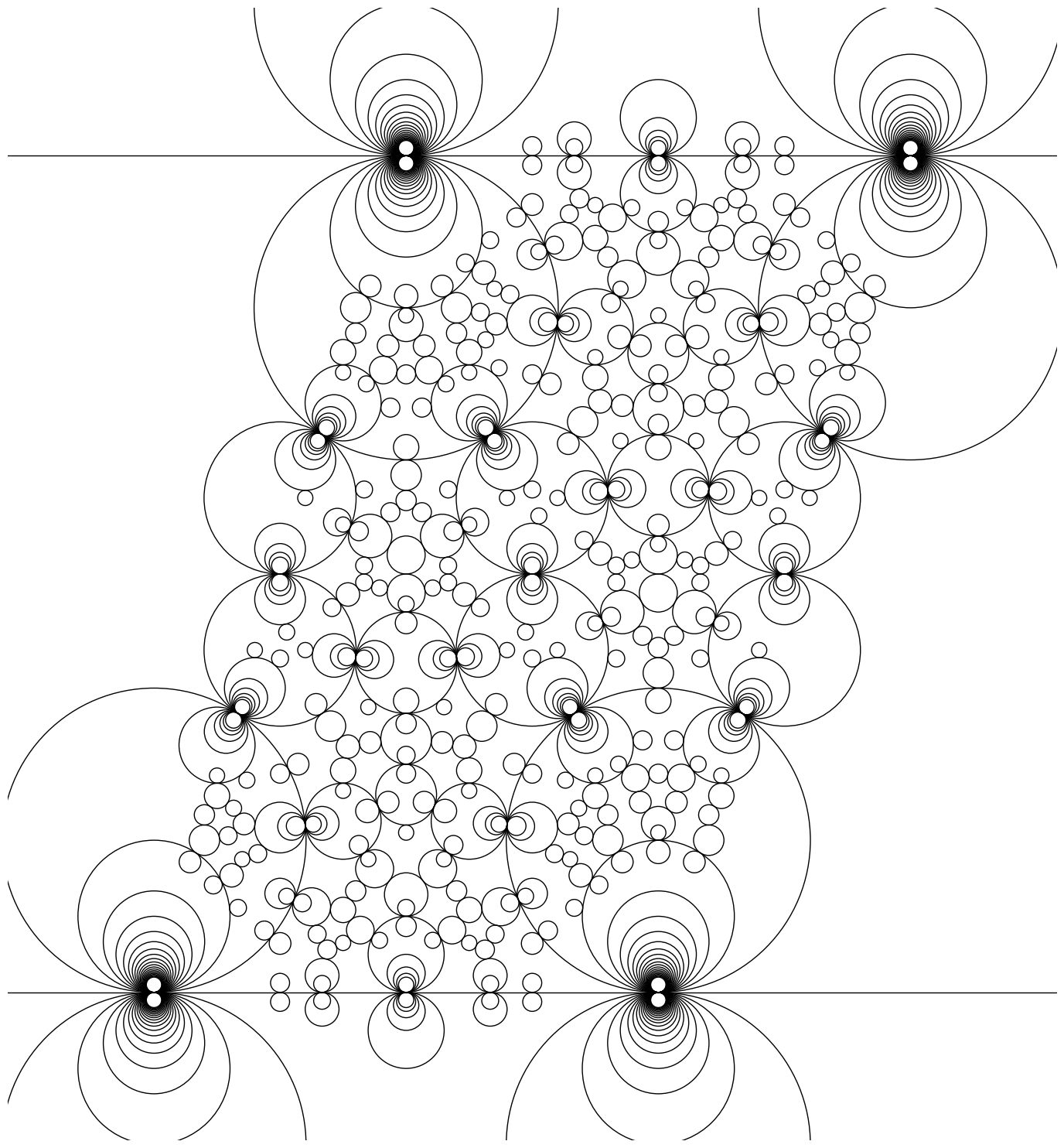} \quad
        \includegraphics[height=2.9in]{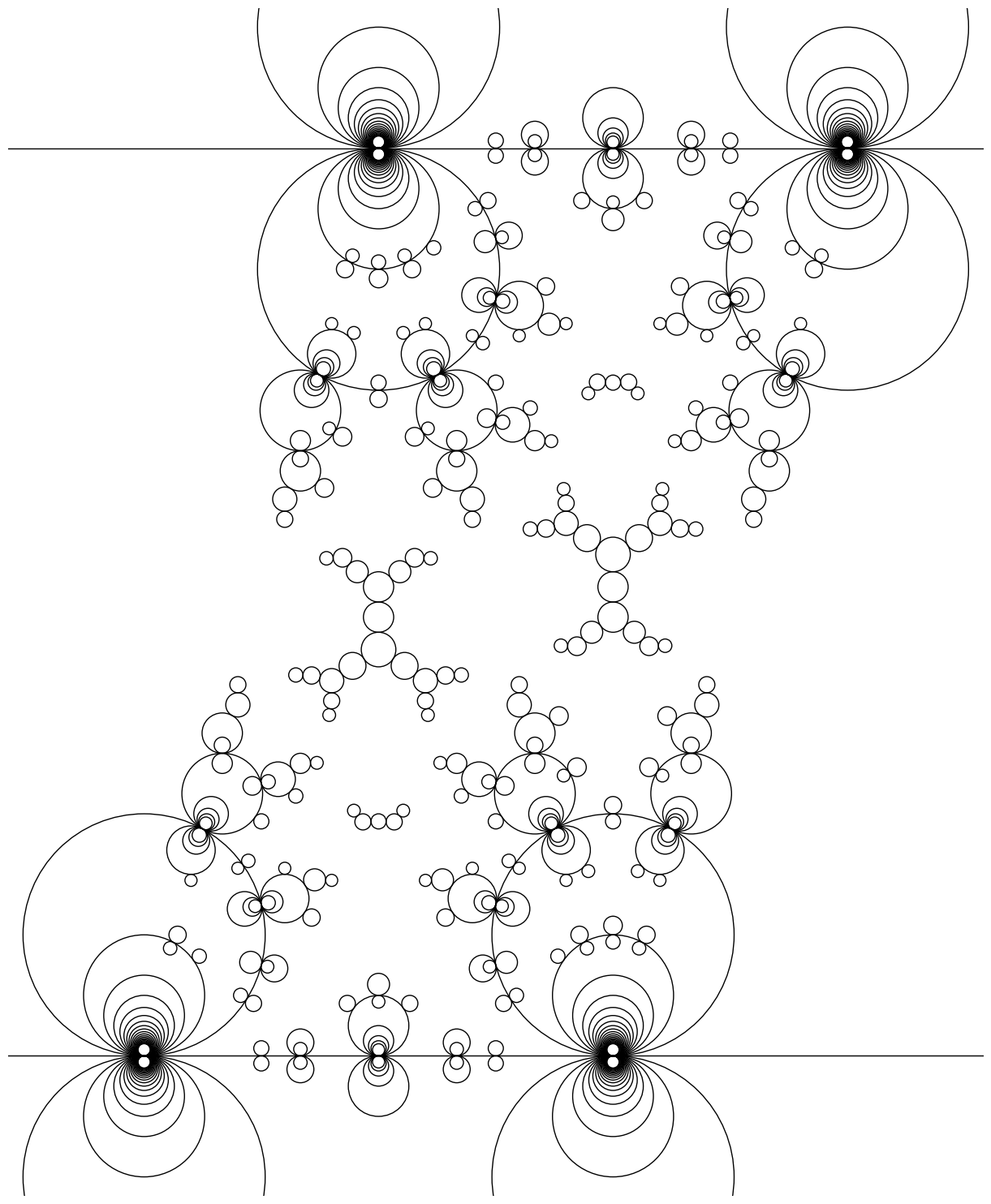} \\
        \caption{Schmidt arrangements $\SK$ of imaginary quadratic fields $K$.  Clockwise from top left:  $\QQ(\sqrt{-2})$, $\QQ(\sqrt{-7})$, $\QQ(\sqrt{-15})$, $\QQ(\sqrt{-11})$.  In each case, the image includes those circles of curvature bounded by $20$ intersecting the fundamental parallelogram of $\OK$ or its boundary.}
\label{fig:exampleSK2}
\end{figure}

Of particular beauty is the orbit of the extended real line $\widehat{\RR} = \RR \cup \{ \infty \}$ under a Bianchi group, which takes the form of an intertwined collection of circles dense in $\widehat{\CC}$, each having curvature (inverse radius) an integer multiple of $\sqrt{-\Delta}$ (where $\Delta$ is the discriminant of $K$).  The circles may intersect only at angles determined by the unit group of $\OK$; only tangently in the case of unit group $\{ \pm 1 \}$ (see Sections \ref{sec:properties} and \ref{sec:intersection} for these and other basic properties).  See Figures \ref{fig:exampleSK3}--\ref{fig:exampleSK2} for examples.  We call the collection a \emph{Schmidt arrangement}, for its appearance, in the cases of discriminant $-3,-4,-7,-8,-11$ in Asmus Schmidt's work on complex continued fractions \cite{\SchmidtFarey,\Schmidt,\SchmidtEleven,\SchmidtEisenstein,\SchmidtTwo}.  

\begin{definition}
Let $K$ be an imaginary quadratic field with ring of integers $\OK$.
The \emph{Schmidt arrangement of $K$}, denoted\footnote{We will use the notation $\Scal_K$ for a subset of $\widehat{\CC}$ and for a collection of circles, without fear of confusion.} $\mathcal{S}_K$, is the orbit of the real line $\widehat{\RR}$ in $\widehat{\CC}$ under the group $\PSL_2(\OK) \le \PSL_2(\CC)$.  The circles in the orbit, i.e. the individual images of $\widehat{\RR}$, are called \emph{$K$-Bianchi circles}.
\end{definition}

The purpose of this paper is to explore the ways in which the arithmetic of $\OK$ shines through in the structure of $\SK$.

In fact, Schmidt's work is the first example of the phenomenon.  The Schmidt arrangement $\Scal_K$ was envisioned by Schmidt as a generalization of the theory of Farey fractions and continued fractions to imaginary quadratic fields.  The relation between Bianchi groups and such complex continued fractions is in direct analogy to the connection between $\PSL_2(\ZZ)$ and real continued fraction expansions, and in this author's opinion is the most natural such generalization she has encountered\footnote{For example, the points of tangency of the collection of circles with curvatures bounded above by $N$ are exactly those elements of $K$ which can be expressed with a denominator bounded in norm by $2N$; this is an analogue of the Farey sequence.}.  

The Schmidt arrangement is also closely related to the action of the Bianchi group as a collection of isometries of the hyperbolic space $\HH^3$.  The extended complex plane $\widehat{\CC}$ is identified with the boundary of $\HH^3$, and a fundamental domain for $\PSL_2(\OK)$ in $\HH^3$ has $h(K)$ cusps on the boundary, where $h(K)$ is the class number of $K$.  
The Schmidt arrangement is the orbit of $\widehat{\RR}$, and its tangency points form the orbit of the cusp corresponding to the trivial class.  Northshield studies the related orbit of cuspidal horoballs \cite{northshield}.

It turns out that the Bianchi groups, among groups $\PSL_2(\OK)$ for all rings of integers $\OK$, have exceptional properties.  For instance, $\PSL_2(\OK)$ contains non-congruence subgroups of finite index when $K=\QQ$ or $K$ is imaginary quadratic (in contrast to all other $K$) \cite{\Lubotzky, \Serre}.  Fine, who, together with Frohman, showed that for $\Delta \neq -3$ the Bianchi groups are non-trivial free products with amalgamation \cite{\FineAmalgam}, points out that the structure of the Bianchi groups is greatly influenced by the number theory of the ring of integers, specifically by whether the ring is Euclidean \cite[\S1.3]{\FineBook}.  Alone among the class of $\PSL_2(\OK)$ for number fields $K$, the Bianchi groups for non-Euclidean imaginary quadratic $\OK$ fail to be generated by elementary matrices \cite{\Unreasonable}.

Schmidt arrangements attracted the author's interest in connection with Apollonian circle packings:  $\SK$ also appears as an `Apollonian superpacking' in \cite{\gttwo}.  Some of the results of this paper generalise results about Apollonian packings to be found in \cite{\ntone, \gttwo}.  This relationship is explored in a companion paper, which is devoted to the study of $K$-Apollonian packings living in $\SK$ and analogues of the Apollonian group \cite{\otherpaper}.  The results of this paper will prove useful in that analysis.

The main results of the present paper form three examples of the arithmetic of $\OK$ appearing in the structure of $\SK$.  

\subsection{First example: curvatures and norms}

Our first example of the arithmetic in $\SK$ is modest:  we illuminate the recursive structure of the circle arrangement in terms of the norm form.   
Fixing one circle in $\SK$, one can describe the curvatures of all tangent circles using the norm form of $\OK$.

\begin{theorem}[Introductory Version of Theorem \ref{thm:quadraticform}]
        Let $K$ be an imaginary quadratic field.  If $C$ and $C'$ are immediately tangent\footnote{See Definition \ref{defn:straddle} for the definition of `immediately.'} $K$-Bianchi circles, then the sum of their curvatures is $\sqrt{-\Delta}N(\beta)$ where $\beta$ is the reduced denominator of the point of tangency.
\end{theorem}

That the circles tangent to one fixed circle in an Apollonian circle packing had curvatures that were values of a translated quadratic form was observed in \cite{\ntone, \Northshield, SarnakLetter} and used as a principle tool in several major results toward conjectures on curvatures in Apollonian packings \cite{\BourgainFuchs,\BourgainKontorovich}.  An excellent exposition is contained in \cite{\FuchsBulletin}.

See Section \ref{sec:intersection} for details.

\subsection{Second example: circles and ideal classes}

The second example provides a deeper explanation for the first:  one can associate certain ideal classes of orders of $\OK$ to the circles in $\SK$. Denote the order of conductor $f$ in $\OK$ by $\Of$.  Write $\Picf$ for the class group of $\Of$, the quotient of its invertible fractional ideals by principal fractional ideals.  There is a homomorphism of groups
\[
        \theta_f: \mathcal{P}ic(\Of) \rightarrow \mathcal{P}ic(\OK)
\]
which is defined by ideal extension: $\theta_f([\mathfrak{a}]) = [\mathfrak{a}\OK]$ for an ideal $\mathfrak{a}$ of $\Of$.

\begin{definition}
        Two $K$-Bianchi circles are \emph{equivalent} (denoted $\sim$) if one can be transformed to the other by $\OK$-translations combined with maps of the form $z \mapsto uz$ for a unit $u$ of $\OK$.
\end{definition}

Then our second exhibit is a natural bijection between the $K$-Bianchi circles modulo equivalence and certain ideal classes.

\begin{theorem}
        \label{thm:mainbij}
        There is a bijection of sets
        \[
                \SK/\sim \;\; 
                \longleftrightarrow \;\; \bigcup_{f \in \ZZ^{>0}} \operatorname{ker}\theta_f
        \]
        in such a way that the circles in bijection with $\operatorname{ker}\theta_f$ have curvature $\sqrt{-\Delta}f$.
\end{theorem}

The relationship between curvatures and conductors allows one to count the circles of given curvature by comparing to the ideal class number of the corresponding order.  In the case of $K=\QQ(i)$, this recovers Theorem 4.2 on Apollonian circle packings of \cite{\ntone} in the context of Theorem 6.3 of \cite{\gttwo}.

It is tantalyzing to ask what `tangency' of ideal classes might represent in the arithmetic of $\OK$; at the moment the author has no particularly compelling answer to this (beyond what is to be found in this paper).

\subsection{Third example: connectedness and Euclideanity}

Our final example concerns the topological structure of $\SK$.

\begin{theorem}[Introductory version of Theorem \ref{thm:connectedEuclidean}]
        Suppose $K$ is a imaginary quadratic field.  Then $\OK$ is Euclidean if and only if $\SK$ is connected.
\end{theorem}

{\bf A note on the figures.} The images in this paper were produced with Sage Mathematics Software \cite{Sage}.

\section{Notations}

Throughout the paper, $K$ is a imaginary quadratic field with discriminant $\Delta<0$ and ring of integers $\OK$.  The ring $\OK$ has an integral basis $1, \tau$, where
\[
        \tau = 
        \left\{\begin{array}{ll}
                 \frac{\sqrt{\Delta}}{2} & \Delta \equiv 0 \pmod 4 \\
               \frac{1 +\sqrt{\Delta}}{2} & \Delta \equiv 1 \pmod 4 \\
        \end{array}
        \right. .
\]
Write $N:K \rightarrow \QQ$ for the norm map $\alpha \mapsto \alpha\overline{\alpha}$.  

\section{Preliminaries on oriented circles in $\widehat{\CC}$}
\label{sec:properties}

Although the main results in the introduction are phrased in terms of (unoriented) $K$-Bianchi circles, it will be convenient to introduce the notion of \emph{orientation}.

\begin{definition}
        An \emph{oriented circle} is a circle together with an \emph{orientation}, which is a direction of travel, specified as either \emph{positive/counterclockwise} or \emph{negative/clockwise}.   An \emph{oriented $K$-Bianchi circle} is an oriented circle whose underlying circle is a $K$-Bianchi circle.  The \emph{interior} of an oriented circle is the area to your left as you travel along the circle according to its orientation.  For lines (circles through $\infty$) besides $\widehat{\RR}$, positive orientation indicates travel in the direction of increasing imaginary part.  For $\widehat{\RR}$, positive orientation indicates travel to the right.
\end{definition}

The action of $\PGL_2(\CC)$ as M\"obius transformations is transitive as an action on circles (since it is transitive on triples of points and three points determine a circle).  It can also be considered to act on oriented circles by taking interiors to interiors.  This action is also transitive (since, for example, $\begin{pmatrix} 0 & 1 \\ 1 & 0 \end{pmatrix}$ reverses the orientation of $\widehat{\RR}$). 

\begin{definition}
        Let us set the convention that $\widehat{\RR}$, when considered an oriented circle, denotes the positively oriented circle (whose interior is the upper half plane).  Any $M \in \PGL_2(\CC)$ has an orientation given by the orientation of $M(\widehat{\RR})$.
\end{definition}
Note that these conventions differ from \cite{\sensual}.

The stabilizer of $\widehat{\RR}$ under the action on circles is $\PGL_2(\RR)$, while the stabilizer under the action on oriented circles is $\PSL_2(\RR)$.

\begin{definition}
        Write $\SK$ for the collection of unoriented $K$-Bianchi circles, and $\OSK$ for the collection of oriented $K$-Bianchi circles.
\end{definition}
There is a $2$-to-$1$ map
        $\OSK \rightarrow \SK$
`forgetting orientation.'  

\begin{proposition}
        Let $K$ be a imaginary quadratic field.  The action of $\PGL_2(\CC)$ on circles described above restricts to a transitive action of $\PSL_2(\OK)$ on $\SK$, with 
        \[
                \operatorname{Stab}(\widehat{\RR}) = \left\{
       \begin{array}{ll}
                \PSL_2(\ZZ) & K \neq \QQ(i) \\
        \left\langle \PSL_2(\ZZ), \begin{pmatrix} 0 & i \\ i & 0 \end{pmatrix} \right\rangle & K = \QQ(i) \\
        \end{array}
        \right. . 
\]  In other words, the circles of $\SK$ are exactly parametrized by cosets $\PSL_2(\OK) / \operatorname{Stab}(\widehat{\RR})$.

Furthermore, the action restricts to a transitive action of
\[
        \widehat{G} =  \left\{
        \begin{array}{ll}
                \PGL_2(\OK) & K \neq \QQ(i) \\
                \PSL_2(\OK) & K = \QQ(i) \\
        \end{array}
        \right.  .
\]
on $\OSK$ where $\operatorname{Stab}(\widehat{\RR}) = \PSL_2(\ZZ)$.  In other words, the oriented circles of $\OSK$ are exactly parametrized by cosets $\widehat{G}/\PSL_2(\ZZ)$.
\end{proposition}

\begin{proof}
Let $\Gamma \subset \PGL_2(\CC)$.  The collection of oriented circles in $\Gamma(\widehat{\RR})$ is parametrized by 
$
        \Gamma / \left(\Gamma \cap \PSL_2(\RR)\right),
$
while the collection of unoriented circles in the orbit is parameterized by
$
        \Gamma / \left(\Gamma \cap \PGL_2(\RR)\right).
$
Taking $\Gamma = \PSL_2(\OK)$ recovers the first statement of the proposition, since $\Gamma$ is transitive on $K$-Bianchi circles, and 
\[
      \PSL_2(\OK) \cap  \PGL_2(\RR)  =
        \left\{ \begin{array}{ll}
                \PSL_2(\ZZ) & K \neq \QQ(i) 
                        \\
                        \left\langle
                \PSL_2(\ZZ), \begin{pmatrix} 0 & i \\ i & 0 \end{pmatrix} \right\rangle & K = \QQ(i) 
                \end{array} \right. .
\]
We now wish to consider $\OSK$.  In the case that $K = \QQ(i)$, the orbit of $\widehat{\RR}$ under $\PSL_2(\OK)$ already obtains both orientations of each $K$-Bianchi circle, since precomposition by $\begin{pmatrix} 0 & i \\ i & 0 \end{pmatrix} \in \PSL_2(\ZZ[i])$ reverses orientation.  Therefore, we need only compute
\[
        \PSL_2(\ZZ[i]) \cap \PSL_2(\RR) 
        = \PSL_2(\ZZ).
\]
However, in the case that $K \neq \QQ(i)$, to obtain both orientations of each $K$-Bianchi circle requires that we consider the orbit of $\widehat{\RR}$ under $\Gamma = \PGL_2(\OK)$.  To see this, note that the negatively oriented real line may not be obtained in $\PSL_2(\OK)$, since $\PSL_2(\OK) \cap \PGL_2(\RR) = \PSL_2(\ZZ)$, no element of which reverse orientation.  On the other hand, since
\[
\Gamma = \PGL_2(\OK) = \PSL_2(\OK) \cup \PSL_2(\OK)\begin{pmatrix} 0 & 1 \\ 1 & 0 \end{pmatrix},        
\]
we see that $\Gamma$ acts transitively on $\OSK$.  Finally, we must compute
\[
        \PGL_2(\OK) \cap \PSL_2(\RR) = \PSL_2(\ZZ).
\]
\end{proof}

\begin{figure}
        \includegraphics[width=4.5in]{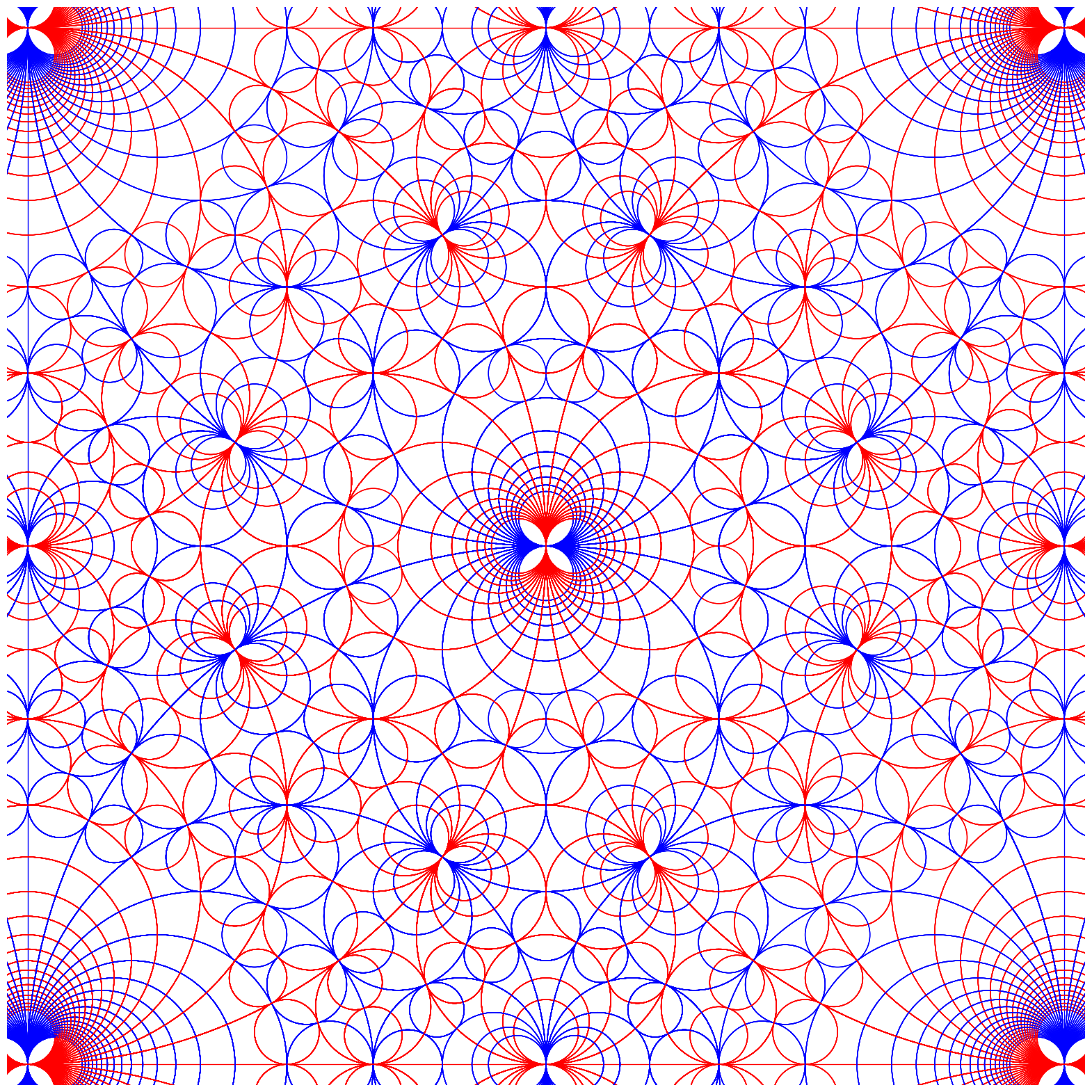}
        \caption{$\OSK$ for $K = \QQ(i)$, showing the region $\{ s + it : 0 \le s,t \le 1\}$.  Red indicates images under $\PSL_2(\OK)$, while blue indicates images under the non-trivial coset of $\PSL_2(\OK)$ in $\PGL_2(\OK)$.}
\label{fig:icosets}
\end{figure}

In the case that $K = \QQ(i)$, the orbit of (unoriented or oriented) circles in $\widehat{\RR}$ under $\PSL_2(\OK)$ is strictly smaller than the orbit under $\PGL_2(\OK)$.  We prefer to study this simpler, smaller orbit, hence the restriction of our definition of $K$-Bianchi circles to $\PSL_2(\OK)$.  The full orbit consists of two copies of the smaller orbit, at right angles (see Figure \ref{fig:icosets}).  In all other $K$, the orbits under $\PGL_2(\OK)$ and $\PSL_2(\OK)$ agree.

\begin{proposition}
        \label{prop:gencirc}
        Let $C$ be an oriented circle in $\widehat{\CC}$ (including those through $\infty$).  Then the circle $C$ can be given uniquely in the form
        \[
                \left\{ X/Y \in \widehat{\CC} : bX\overline{X} - a Y \overline{X} - \overline{a} X \overline{Y} + b' Y \overline{Y} = 0 \right\}.
        \]
        where the following hold:
        \begin{enumerate}
                \item \label{item1}$b,b' \in \RR$, $a\in \CC$,
                \item we have
                        \begin{equation}
                \label{eqn:cocurv}
                b'b = a\overline{a}-1,
        \end{equation}
\item\label{item3} $b$ has sign equal to the orientation of $C$, and;
\item\label{item4} if $b=0$, in which case $C$ must be a line, then $a$, as a vector, is a unit vector pointing from exterior to interior orthogonal to $C$.
        \end{enumerate}
        Furthermore, if $C'$ is the image of $C$ under $z \mapsto 1/z$, and $C$ has parameters $(b,b',a)$ according to the requirements above, then $C'$ has parameters $(b',b,\overline{a})$ according to the requirements above.
        Finally, \begin{enumerate}
                \item $b=0$ if and only if $\infty \in C$,
                \item $b'=0$ if and only if $0 \in C$,
                \item if $\infty \notin C$, then $C$ has radius $1/|b|$,
                \item if $\infty \notin C$, then $C$ has centre $a/b$.
        \end{enumerate}
\end{proposition}

\begin{proof}
        When $b \neq 0$, using $b' := \frac{a\overline{a} -1}{b}$, the equation given is exactly the equation
        \[
                \left| \frac{a}{b} - \frac{X}{Y} \right| = \frac{1}{|b|}.
        \]
        This equation uniquely determines an oriented circle not passing through $\infty$.  The conditions on $a,b,b'$ in the statement ensure that no two equations give the same oriented circle.

        Setting $b=0$, but assuming $b' \neq 0$, one obtains an equation of the form
        \[
                \operatorname{Re}\left( \overline{a} \frac{X}{Y} \right) = b'
        \]
        under the condition $a\overline{a}=1$.  This uniquely determines an oriented circle through $\infty$ but not passing through $0$, and no two such equations give the same oriented circle.

        Finally, setting $b=b'=0$, one obtains
        \[
                \operatorname{Re}\left( \overline{a} \frac{X}{Y} \right) = 0 
        \]
        under the condition that $a\overline{a}=1$ and $a$ satisfies the argument condition given in the statement.  This uniquely determines a circle through $0$ and $\infty$, and no two such equations give the same oriented circle.

        The other statements of the proposition are all immediate consequences of these observations, or straightforward verifications.
\end{proof}

\begin{definition}
        \label{def:curv}
        For any circle $C$ in $\widehat{\CC}$, expressing it as in Proposition \ref{prop:gencirc}, we call $b$ the \emph{curvature} (elsewhere sometimes called a \emph{bend}), $b'$ the \emph{co-curvature}, and $a$ the \emph{curvature-centre}.
\end{definition}

The following proposition appears for the case of $\ZZ[i]$ in \cite{\sensual}.

\begin{proposition}
\label{prop:curvature}
Consider an oriented circle expressed as the image of $\widehat{\RR}$ under a transformation of the form
\[
M =
      \begin{pmatrix} \alpha & \gamma \\ \beta & \delta \end{pmatrix}, \quad \alpha, \beta, \gamma, \delta \in \CC, \;\; |\alpha\delta - \beta\gamma|=1.
\]
The curvature of the circle is given by
\[
        i(  \beta \overline{\delta} - \overline{\beta}\delta),
\]
the co-curvature of the circle is given by
\[
        i( \alpha \overline{\gamma} - \overline{\alpha}\gamma ),
\]
and the curvature-centre is given by
\[
        i(  \alpha \overline{\delta} - \gamma \overline{\beta} ).
\]
\end{proposition}

\begin{proof}
        To prove the formul{\ae}, it suffices to verify that the three points
        \begin{equation}
                \label{eqn:z}
                \frac{X}{Y} = \frac{\alpha}{\beta}, \frac{\gamma}{\delta}, \frac{\alpha+\gamma}{\beta+\delta}
   \end{equation}
   lie on the circle
   \begin{equation}
           \label{eqn:circ}
                 bX\overline{X} - a Y \overline{X} - \overline{a} X \overline{Y} + b' Y \overline{Y} = 0 
   \end{equation}
   and that the quantities satisfy the extra requirements of Proposition \ref{prop:gencirc}.  The verifications of \eqref{eqn:circ} and \eqref{eqn:cocurv} are straightforward and use the fact that $N(\alpha\delta - \gamma\beta)=1$.  Requirement \eqref{item1} of Proposition \ref{prop:gencirc} is also immediately verified.

   We must demonstrate that orientation matches the sign of the formula for curvature.  The map
   \[
   \phi: \PSL_2(\CC) \rightarrow \RR, \quad \begin{pmatrix} \alpha & \gamma \\ \beta & \delta \end{pmatrix} \mapsto i(\beta\overline{\delta}-\overline{\beta}\delta)
   \]
   is surjective and continuous.  Thus, removing $\phi^{-1}(0)$ necessarily leaves two path connected components of $\PGL_2(\CC)$ upon which the sign of $\phi(M)$ is constant -- one with positive sign and one with negative sign.  
   
   The vanishing of $i(\beta\overline{\delta} - \overline{\beta}\delta)$ is equivalent to $\beta = \lambda\delta$ for some $\lambda \in \RR$; this, in turn is equivalent to the existence of a non-trivial vanishing $\RR$-linear combination of $\beta$ and $\delta$; this, in turn, is exactly the statement that $M$ takes some element of $\widehat{\RR}$ to $\infty$.  Hence straight lines are exactly those circles for which $i(\beta\overline{\delta} - \overline{\beta}\delta)$ vanishes.

   Therefore, $\phi^{-1}(0)$ consists of all matrices taking $\widehat{\RR}$ to a circle through $\infty$ (i.e. a line).   It follows that one component of $\PGL_2(\CC) \setminus \phi^{-1}(0)$ consists of those $M$ taking $\widehat{\RR}$ to circles of positive orientation, and the other those taking $\widehat{\RR}$ to circles of negative orientation.  For, orientation of circles cannot reverse under continuous deformation without passing through $\infty$.  
  
  It remains to verify that the orientation matches the sign of the claimed formula for curvature for at least one example; any matrix $M$ from one of the two components will suffice, e.g. $\begin{pmatrix} 1 & 0 \\ i & 1 \end{pmatrix}$.  The associated transformation maps $\widehat{\RR}$ to a circle of radius $1/2$ centred on $-i/2$, taking $0 \mapsto 0$, $1 \mapsto (1-i)/2$, $\infty \mapsto -i$.  Hence it has negative orientation and $i(\beta\overline{\delta}-\overline{\beta}\delta) = -2$.

   Finally, we must address the special case of circles passing through $0$ and $\infty$.  Such a circle can be expressed in a particularly simple form as $M(\widehat{\RR})$ for
   \[
 M =           \begin{pmatrix}
         u & 0 \\
 0 & 1 \end{pmatrix} \in \PGL_2(\OK), u \in \OK^*.
 \]
 Then the formula for curvature-centre gives $a = iu$ and requirement \eqref{item4} of Proposition \ref{prop:gencirc} regarding this special case is verified.
\end{proof}

Note that the curvature and co-curvature of a $K$-Bianchi circle are the real parts of elements of $2i\OK$, so they are integer multiples of $\sqrt{-\Delta}$; the integer alone will be referred to as the \emph{reduced curvature} or \emph{reduced co-curvature}, respectively.

There is a natural way to view circles as points in Minkowski space \cite{\gtone, \Kocik}.  Let $\MM$ be Minkowski space, that is, the vector space $\RR^4$ endowed with inner product 
\[
        \langle (x_1, x_2, x_3, x_4) , (y_1, y_2, y_3, y_4) \rangle = - x_1y_2 - x_2y_1 + x_3y_3 + x_4y_4.
\]
Define the Pedoe map
\[
        \pi: \{ \mbox{oriented circles} \} \rightarrow \MM
        \]
        by letting $\pi(C) = (b,b',x,y)$ where $b$ and $b'$ are the curvature and co-curvature of $C$, and $x+iy = z$ is the curvature-centre.  By Proposition 4.1, $\langle \pi(C), \pi(C) \rangle = 1$, so that the image of the Pedoe map lies in the hypersurface $Q_M=1$ where $Q_M$ is the quadratic form of $\MM$.  Then the inner product gives an inner product on circles carrying geometric information.

\begin{figure}
        \includegraphics[height=1in]{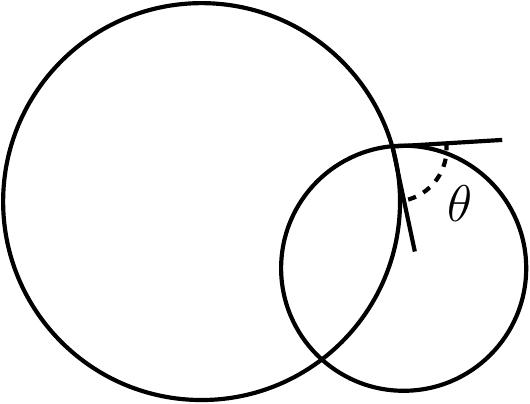} 
        \caption{The angle of intersection of two circles.}
\label{fig:anglecircles}
\end{figure}

\begin{proposition}[Proposition 2.4 of \cite{\Kocik}]
        \label{prop:pedoeproduct}
        Let $v_i = \pi(C_i)$ for two circles $C_1, C_2$ which are not disjoint.  Then $\langle v_1, v_2 \rangle = \cos \theta$, where $\theta$ is the angle between the two circles as in Figure \ref{fig:anglecircles}.  In particular,
                \begin{enumerate}
                        \item $\langle v_1, v_2 \rangle = -1$ if and only if the circles are tangent externally \includegraphics[width=0.5in]{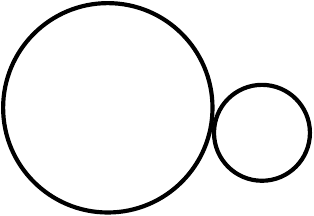}
                        \item $\langle v_1, v_2 \rangle = 1$ if and only if the circles are tangent internally \includegraphics[width=0.38in]{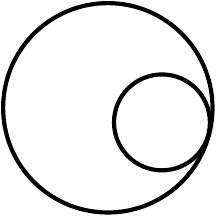}
                        \item $\langle v_1, v_2 \rangle = 0$ if and only if the circles are mutually orthogonal \includegraphics[width=0.47in]{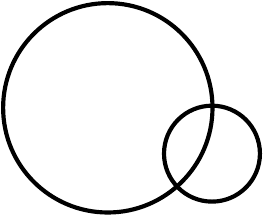}
                \end{enumerate}
\end{proposition}

Note that the Pedoe product of two $K$-Bianchi circles lies in $\frac{1}{2}\ZZ$ (see Proposition \ref{prop:curvature}).  This indicates that the angles of intersection of $K$-Bianchi circles will be restricted, as we will see in the next section.

\section{The intersection of $K$-Bianchi circles}
\label{sec:intersection}

The delicate arrangement of circles in $\Scal_K$ is controlled by the number theory of $K$.  In this section, we examine some of their basic properties:
\begin{enumerate}
        \item We show that oriented $K$-Bianchi circles intersect only at $K$-points and only in angles prescribed by the unit group of $\OK$ (usually only tangently).  
        \item We describe exactly the circles passing through any point $z \in K$.
        \item We give some sufficient and necessary conditions for circles of various curvatures to be tangent or pass through a given point.  
\end{enumerate}

\begin{proposition}
        \label{prop:Kintersect}
        $K$-Bianchi circles intersect only at $K$-points.
\end{proposition}

\begin{proof}
        Without loss of generality, we may consider one of any two intersecting circles to be $\widehat{\RR}$.  Let the other circle $C$ be $M(\widehat{\RR})$ for some $M \in \PGL_2(\OK)$.  Suppose $s$ is an intersection point of $C$ and $\widehat{\RR}$.  Then there exists some $t$ of $\widehat{\RR}$ such that $M(t) = s \in \widehat{\RR}$.  In other words, $\overline{M}(t) = \overline{M}(\overline{t}) = M(t)$.  Without loss of generality, we suppose that $t \neq \infty$ (for example, precompose $M$ with $z \mapsto 1/z$).  Suppose
        \[
                M = \begin{pmatrix} \alpha & \gamma \\ \beta & \delta \end{pmatrix}.
        \]
        Then $M(t) = \overline{M}(t)$, scaled by $\frac{1}{2i}$, becomes a quadratic equation
        \begin{equation}
                \label{eqn:quadkintersect}
                \operatorname{Im}(\alpha \overline{\beta}) t^2 + \operatorname{Im}(\gamma \overline{\beta} + \alpha \overline{\delta}) t + \operatorname{Im}(\gamma \overline{\delta}) = 0.
        \end{equation}
        Let the transpose $M^t \in \PGL_2(\OK)$ have curvature-centre $z = x+iy$, curvature $b$ and co-curvature $b'$.  Then the coefficients of the quadratic \eqref{eqn:quadkintersect} are, respectively, $b'$, $x$, and $b$; these are all in $\frac12\sqrt{-\Delta}\ZZ$ (by Proposition \ref{prop:curvature}).  The discriminant of \eqref{eqn:quadkintersect} (using \eqref{eqn:cocurv}) is
        \[
x^2 - 4bb' = 4 - 3(x^2+y^2) - y^2
\]
Since $t \in {\RR}$, this discriminant is non-negative.  As $z=x+iy \in i\OK$, we have $x^2 + y^2 \in \ZZ$.  Thus the possible non-negative values of the discriminant are:
\[
        4, \frac{9}{4}, 1, \frac{1}{4}, 0,
\]
and in any case, $t$ must be a rational root, hence $s$ is rational.
We return to the general case by applying a M\"obius transformation with coefficients in $K$ and find that the intersection points lie in $K$.
\end{proof}

In the proof, one may remark that in fact the discriminant of \eqref{eqn:quadkintersect} lies in $\frac{\Delta}{4}\ZZ$.  Hence, for large $-\Delta$, the discriminant is exactly $0$ and the circles are necessarily tangent.  The following proposition strengthens this observation.

\begin{proposition}
        \label{prop:angles}
        Two $K$-Bianchi circles may intersect at an angle of $\theta$ only if $e^{i\theta} \in \OK$.  In particular, if the unit group of $\OK$ is $\{ \pm 1\}$, then circles may only be disjoint or tangent.
\end{proposition}

\begin{proof}
        Consider the circles passing through $0$.  These are all images of $\widehat{\RR}$ under M\"obius transformations of the form
\[
        \begin{pmatrix} 0 & \gamma \\
                        1 & \delta
\end{pmatrix} \in \PGL_2(\OK),
\]
where $\gamma$ is a unit in $\OK$.  
These have curvature-centre $i\gamma$.
Therefore the centre lies on the lines through $0$ which pass through the units of $\OK$.  
Hence the circles through $0$ may meet at angles $\theta$ for which $e^{i\theta} \in \OK$. 

Since the intersection points of any two $K$-Bianchi circles lie in $K$ (by Proposition \ref{prop:Kintersect}), and M\"obius transformations preserve angles, we may transport the point of intersection to the origin to prove the general statement.
\end{proof}

For the case $K = \QQ(i)$, the proof is actually slightly stronger than the statement:  it applies to any circles which are images of $\widehat{\RR}$ under $\PGL_2(\OK)$.

We say that $\alpha, \beta \in \OK$ are coprime if the ideals they generate are coprime, i.e. $(\alpha)+(\beta)=(1)$.

\begin{proposition}
        \label{prop:tangentfamilies}
        Let $\alpha/\beta \in K$ be such that $\alpha$ and $\beta$ are coprime.  Suppose $|\OK^*| = n$.  Then the collection of oriented $K$-Bianchi circles passing through $\alpha/\beta$ is a union of $n$ generically different $\ZZ$-families, one for each $u \in \OK^*$.  The family associated to $u$ consists of the images of $\widehat{\RR}$ under the transformations
\[
        \begin{pmatrix}
                \alpha & u \gamma + k \tau  \alpha\\
                 \beta & u \delta + k \tau \beta 
        \end{pmatrix}, \quad k \in \ZZ,
\]
where $\gamma$, $\delta$ is a particular solution to $\alpha \delta - \beta\gamma = 1$.  Furthermore,
\begin{enumerate}
        \item The curvatures of the circles in one family form an equivalence class modulo $\sqrt{-\Delta}N(\beta)$.  
        \item The centres of the circles in a given family lie on a single line through $\alpha/\beta$.  
        \item The family given by unit $u$ contains the same circles as the family given by $-u$, but with opposite orientations.
\end{enumerate}
\end{proposition}

\begin{proof}
The circles passing through $\alpha/\beta$ are those with representative matrices of the form
\[
        \begin{pmatrix}
                \alpha & \gamma \\
                 \beta & \delta 
        \end{pmatrix}
\]
where $\alpha \delta - \gamma\beta = u$ is a unit in $\OK$ (and $\alpha, \beta, \gamma, \delta \in \OK$).  There is a particular solution, call it $\gamma, \delta$, for the unit $u=1$, and homogeneous solutions $\gamma + \eta \alpha$, $\delta + \eta\beta$ for $\eta \in \OK$.  Since the circle does not depend on $\eta$ up to addition of a rational integer, we obtain a family of circles
\[
        \begin{pmatrix}
                \alpha & \gamma + k \tau  \alpha\\
                 \beta & \delta + k \tau \beta 
        \end{pmatrix}, \quad k \in \ZZ.
\]
The curvatures of this family are
$i \left( \beta \overline{\delta} - \overline{\beta}\delta \right) + k(\sqrt{-\Delta}) N(\beta)$
and therefore they are all distinct circles.  Their curvature-centres are
$i \left( \gamma \overline{\beta} - \alpha\overline{\delta}\right) + k(\sqrt{-\Delta})\alpha\overline{\beta}$.
The difference between two centres (say, $k=0$ and $k=k_0$) has the form
        $\frac{ k_0 X }{Y(k_0)}$, 
where $X$ doesn't depend on $k_0$, and $Y(k_0)$ is real.  Therefore, the circles in one family all lie with their centres on a single line.  Since their curvatures approach $\infty$, this line must pass through $\alpha/\beta$.  
To prove the analogous results for the general family associated to $u$, replace $\gamma$, $\delta$ with $u\gamma$ and $u\delta$.
To justify the final assertion of the Proposition, note that changing the sign of $u$ changes the signs of the curvature and co-curvature but does not change the collection of centres.
\end{proof}

Combining the previous two propositions, we learn that the the angles between the lines of centres for any two families are given by some $\theta$ where $e^{i\theta}$ is a unit of $\OK$.

\begin{proposition}
        \label{prop:onlytangent}
        If $K \neq \QQ(\sqrt{-3})$ then the circles of $\SK$ may only intersect tangently.
\end{proposition}
\begin{proof}
        For $K \neq \QQ(i)$, this is immediate from Proposition \ref{prop:angles}.  For $K = \QQ(i)$, oriented circles given by $\PGL_2(\OK)$ may intersect orthogonally.  However, from Proposition \ref{prop:tangentfamilies}, two orthogonal families at a point have determinants $\{\pm 1\}$ and $\{ \pm i \}$.  In an equivalence class of $\PGL_2(\OK)$, the determinants are either $\{ \pm 1 \}$ or $\{ \pm i\}$; the subgroup $\PSL_2(\OK)$ consists of those having the former shape.  Therefore, by restricting to circles in $\SK$, one considers only $M(\widehat{\RR})$ for $M \in \PGL_2(\OK)$ of determinant $\pm 1$.  At a single point, we obtain only two of the four families described in Proposition \ref{prop:tangentfamilies} and these all intersect tangently.
\end{proof}

For $K = \QQ(\sqrt{-3})$, the circles of $\SK$ do intersect at angles of $\pi/3$ and $2\pi/3$, for which reason this field presents extra difficulties.

The following definition will be key to describing the tangency structure of $\SK$.

\begin{definition}\label{defn:straddle}
        A collection $\Pcal \subset \SK$ of circles \emph{straddles} a circle $C$ if $\Pcal$ includes a circle $C' \neq C$ disjoint from the interior of $C$ and a circle $C'' \neq C$ disjoint from the exterior of $C$.  
        Two oriented $K$-Bianchi circles are \emph{immediately tangent} if they are externally tangent in such a way that the pair straddles no circle of $\SK$.
\end{definition}

If $K \neq \QQ(\sqrt{-3})$, then $\Pcal$ straddles $C$ exactly if $\Pcal$ intersects both the interior and exterior of $C$.   This is because the restriction that intersections only occur as tangencies (Proposition \ref{prop:onlytangent}) implies that no circle of $\Pcal$ may intersect interior and exterior of $C$ simultaneously.

\begin{proposition}
        \label{prop:immtang}
        Let $C \in \OSK$ be an oriented $K$-Bianchi circle with $K$-rational point $x$.  Then there exists 
\[
        M_C = \begin{pmatrix} \alpha & \gamma \\ \beta & \delta \end{pmatrix} \in \PGL_2(\OK)
\]
such that $C = M_C(\widehat{\RR})$ and $x = \alpha/\beta$.  Furthermore, there exists exactly one oriented $K$-Bianchi circle $C' \in \OSK$ immediately tangent to $C$ at $x$, given by $C' = M_{C'}(\widehat{\RR})$ where
\[
        M_{C'} = 
         \begin{pmatrix} \alpha & -\gamma + \tau \alpha \\ \beta & -\delta + \tau\beta \end{pmatrix}
        = M_C \begin{pmatrix} 1 & \tau \\ 0 & -1 \end{pmatrix} \in \PGL_2(\OK).
\]
\end{proposition}

In the case of $K = \QQ(i)$, these matrices may be taken in $\PSL_2(\OK)$.

\begin{proof}
        Let $M$ be some matrix such that $M(\widehat{\RR}) = C$.  Then the $K$-rational points of $C$ are given by $\alpha/\beta$ where $(\alpha\; \beta)^T$ is a $\ZZ$-linear combination of the column vectors of $M$.  Hence, some $N \in \PSL_2(\ZZ)$ will give the desired $M_C = MN$.
    
        Proposition \ref{prop:tangentfamilies} lists the circles tangent to $C$ at $x$.  There are two families of circles tangent to $C$:  those whose centres lie on the line passing through $x$ and the centre of $C$ (given by $u= \pm 1$ in the language of Proposition \ref{prop:tangentfamilies}).  The circle $C'$ must be a member of one of these two families.  Circles in one family are orientated such that they pass through $x$ in the same direction.  Therefore, the circles must belong to distinct families in order to be externally tangent.  The underlying unoriented family of each of the two families is the same, and $C$ and $C'$ must be `consecutive' in this unoriented family.  That is, $C'$ is given by $k = \pm 1$, $u=-1$ in the language of Proposition \ref{prop:tangentfamilies}.  The sign of $k$ is determined by the requirement that the circles have disjoint interiors (using the wrong sign will mean disjoint exteriors).
\end{proof}

\begin{theorem}
        \label{thm:quadraticform}
        Fix an oriented $K$-Bianchi circle $C \in \OSK$ of reduced curvature $b$, given by $M(\widehat{\RR})$ for
        \[
                M = \begin{pmatrix} \alpha & \gamma \\ \beta & \delta \end{pmatrix}.
        \]
        Let
        \[
                \Lambda = \beta \ZZ + \delta \ZZ \subset \OK.
        \]
        Then the reduced curvatures of the collection of oriented $K$-Bianchi circles immediately tangent to $C$ are exactly the primitive values of the translated integral binary quadratic form
        \[
N(x) - b, \quad x \in \Lambda.
\]
\end{theorem}

\begin{proof}
        Immediate from Proposition \ref{prop:immtang} and the fact that the $K$-rational points on $C$ have denominators chosen from $\Lambda$.
\end{proof}

This result, in the case $K=\QQ(i)$, in the form of a statement about curvatures in Apollonian circle packings, is due to Sarnak \cite{SarnakLetter}.

\section{The bijection between ideal classes and $K$-Bianchi circles}

Tacit in the results of the last section is a study of rank two $\ZZ$-sublattices of $\OK$; this is exactly the study of ideals of orders $\Of$ of $K$, as we will presently explain.

We begin with some background on lattices.  If one lattice is contained in another, $\Lambda \subseteq \Lambda'$, then the index of $\Lambda$ in $\Lambda'$, denoted $[\Lambda': \Lambda]$, is defined to be the size of the quotient $\Lambda'/\Lambda$.  In this paper, we consider rank two $\ZZ$-lattices in $\OK$, and the index in $\OK$ will be called the \emph{covolume}.  If a basis for $\Lambda$ is given by $a+b\omega$, $c+d\omega$, then its covolume is $|ad-bc|$.

In the case that $\Lambda$ is equal to an integral ideal of $\OK$, the covolume is equal to the norm of the ideal.  If it is an integral ideal of an order $\Of$ of conductor $f$, then its covolume is $f$ times the norm.  A fractional ideal class of some $\Of$ is a class of lattices up to homothety (defined as scaling by an element of $\OK$).

The order of any rank two $\ZZ$-lattice $\Lambda \in \OK$ is defined as
\[
   \operatorname{Ord}(\Lambda) =        \{ \alpha \in \OK : \alpha\Lambda \subseteq \Lambda \}.
\]
The lattice $\Lambda$ is always a representative lattice for some fractional ideal class of $\operatorname{Ord}(\Lambda)$, which is an order of $K$.  Homothety preserves the order of a lattice.

\begin{proposition}
        \label{proposition:primeval}
        Let $\Lambda$ be a rank two $\ZZ$-lattice contained in $\mathcal{O}_K$.  The following are equivalent:
        \begin{enumerate}
                \item \label{cond1} $\Lambda$ has a basis of coprime elements.
                \item \label{cond2} $\Lambda$ has conductor equal to its covolume.
        \end{enumerate}
        Furthermore, an ideal class $[\mathfrak{a}]$ of $\Of$ lies in $\operatorname{ker} \theta_f$ if and only if it has a representative lattice $\Lambda$ satisfying the conditions above.  In this case, such a $\Lambda$ is unique up to homothety by a unit of $\mathcal{O}_K$.  Further, if $\beta, \delta$ form a basis for $\Lambda$, then the covolume is equal to $\frac{1}{\sqrt{-\Delta}}|\overline{\beta}\delta - \beta\overline{\delta}|$.
\end{proposition}

We will call such a lattice \emph{primeval}.

\begin{proof}

        First, we show that $[\Lambda] \subseteq \operatorname{ker}\theta_f$ for some $f$ if and only if it has a representative satisfying condition \eqref{cond1} of the statement.  Suppose \eqref{cond1}, and that $[\Lambda]$ is an ideal class of $\Of$.  Then the $\OK$-span of $\Lambda$ is $\OK$, i.e. $\theta_f([\Lambda]) = [\OK]$.  Conversely, if $\theta([\Lambda]) = [\OK]$, then some homothety $\Lambda'$ of $\Lambda$ generates $\OK$, i.e., satisfies \eqref{cond1}.  The uniqueness statement follows from the observation that two lattices in the same fractional ideal class of $\Of$, both lying in $\mathcal{O}_K$ and having covolume $f$, are necessarily related by homothety by a unit of $\mathcal{O}_K$. 

        Now suppose $\Lambda$ has a basis $\beta = a+b\omega, \delta = c+d\omega$ of coprime elements.   From the basis, the covolume of $\Lambda$ is $|ad-bc|$.  It is easy to check that
        \begin{equation}
                \label{eqn:imbd}
                                            \beta\overline{\delta} - \overline{\beta}\delta
                                            = \sqrt{\Delta}(bc-ad),
                                    \end{equation}
                                    from which the last sentence of the Proposition statement follows.
        
        It remains to compute the conductor of $\Lambda$ in terms of its covolume.
        Let $\mathcal{O}_f$ be the order of $\Lambda$, having conductor $f$.   The inverse of the class of $\Lambda$ is its conjugate $\overline{\Lambda}$ in the ideal classes of $\mathcal{O}_f$.  In other words, the product lattice is homothetic to the order itself:
        \[
                \Lambda \overline{\Lambda} = \langle \lambda \overline{\mu} : \lambda, \mu \in \Lambda \rangle = \nu \mathcal{O}_f,
        \]
        for some $\nu \in \mathcal{O}_K$.  
                                   The order $\mathcal{O}_f$ can be characterized as those elements of $\OK$ having imaginary part an integer multiple of $f\frac{\sqrt{\Delta}}{2}$.
                                   Now,
\[
                \Lambda \overline{\Lambda}  = \left\langle N(\beta), N(\delta), \beta\overline{\delta}, \overline{\beta}\delta \right\rangle.
\]
These generators have imaginary parts $0$ and $\pm [\OK:\Lambda]\frac{\sqrt{\Delta}}{2}$ by \eqref{eqn:imbd}.  Furthermore, as the rational integers
        \[
                N(\beta), N(\delta), 
                                            \beta\overline{\delta} + \overline{\beta}\delta
                                    \]
                                    have no common factor above $1$, we find $1 \in \Lambda\overline{\Lambda}$, so $\nu = \pm 1$. Consequently, 
$f = [\OK:\Lambda]$. 
\end{proof}

The maps which witness the bijection of Theorem \ref{thm:mainbij} are quite simply defined.  Denote them
        \[
                I: \SK/\sim\; \longrightarrow\; \bigcup_{f \in \ZZ^{>0}} \operatorname{ker}\theta_f
                \quad\quad\mbox{and}\quad\quad
                S: \bigcup_{f \in \ZZ^{>0}} \operatorname{ker}\theta_f
                \;\longrightarrow \;\SK/\sim.
        \]

        Given $C \in \SK$, let $M \in \PSL_2(\OK)$ take $\widehat{\RR}$ to $C$.  Write
        \[
                M = \begin{pmatrix} \alpha & \gamma \\ \beta & \delta \end{pmatrix}.
        \]
        Then $I(C)$ is the ideal class of the lattice $\ZZ$-generated by $\beta$ and $\delta$.

        For the reverse map, one chooses a primeval representative lattice of the ideal class $\mathfrak{a}$, and selects a $\ZZ$-basis $\beta, \delta$ of coprime integers.  Solving for $\alpha, \gamma \in \OK$ such that
        \[
                M = \begin{pmatrix} \alpha & \gamma \\ \beta & \delta \end{pmatrix}.
        \]
        is an element of $\PSL_2(\OK)$, we set $S(\mathfrak{a}) = M(\widehat{\RR})$.

\begin{proof}[Proof of Theorem \ref{thm:mainbij}]
The matrices
\[
        M = \begin{pmatrix} \alpha & \gamma \\ \beta & \delta \end{pmatrix},
\;        \mbox{and} \;
        M' = \begin{pmatrix} \alpha' & \gamma' \\ \beta' & \delta' \end{pmatrix} \in \PGL_2(\OK)
\]
map $\widehat{\RR}$ to equivalent unoriented circles if and only if $M' = TMS$ where $S \in \PGL_2(\ZZ)$, and $T = \tiny \begin{pmatrix} \epsilon & \eta \\ 0 & 1 \end{pmatrix}\normalsize$ (a translation by $\eta \in \OK$ followed by a multiplication by a unit $\epsilon \in \OK$).  This occurs if and only if $\beta, \delta$ and $\beta', \delta'$ are bases, each consisting of a pair of coprime elements, for the same primeval lattice, up to homothety by a unit.   The primeval lattices form a system of representatives for $\cup_{f\in\ZZ^{>0}} \operatorname{ker}\theta_f$ by Proposition \ref{proposition:primeval}.  These observations suffice to show that the maps $I$ and $S$ are well-defined and inverses.  By Proposition \ref{prop:curvature}, the oriented curvature of the circle $M(\widehat{\RR})$ is $i({\beta}\overline{\delta} - \overline{\beta}{\delta})$, while by Proposition \ref{proposition:primeval}, the conductor of the lattice is $(-\Delta)^{-\frac{1}{2}}|\overline{\beta}\delta - \beta\overline{\delta}|$.
\end{proof}

\section{Counting $K$-Bianchi circles}

Let $h_K$ denote the class number of $K$.  Let $h_f$ denote the class number of an order $\Of$ of conductor $f$ inside $K$ (so $h_1 = h_K$).  Define
\[
        \mathcal{U}_f = \OK^*/\Of^*.
\]
Let $u = |\mathcal{U}_f|$.  This group is trivial except for $\Delta = -3, -4$.  We have
        \[
                u = \left\{ \begin{array}{ll}
                        2 & K=\QQ(\sqrt{-1}) \\
                        3 & K=\QQ(\sqrt{-3}) \\
                        1 & \mbox{otherwise} 
                        \end{array} \right. .
        \]

        There is an exact sequence (see, for example, \cite[{\S I.12}]{\Neukirch}):
\begin{equation}
        \label{eqn:exact}
        \xymatrix{
                1 \ar[r] & \mathcal{U}_f \ar[r] & (\OK/{f})^*/(\ZZ/{f})^* \ar[r] & \mathcal{P}ic(\Of) \ar[r]^{\theta_f} & \mathcal{P}ic(\OK) \ar[r] & 1
}
\end{equation}

Consequently, the following standard result gives the size $h_f$ of $\mathcal{P}ic(\Of)$.

\begin{theorem}[{See, for example, \cite[{\S XIII.2}]{\CohnBook}, \cite[Theorem 7.24]{\Cox} or \cite[{\S I.12}]{\Neukirch}}]
        If $f>1$,
        \[
                h_f = \frac{h_K}{u} f \prod_{\substack{p \mid f \\ p \mbox{\tiny\; prime}}} \left(1 - \frac{1}{p}\left(\frac{\Delta}{p}\right) \right).
        \]
\end{theorem}

As a corollary to Theorem \ref{thm:mainbij}, we obtain

\begin{corollary}[Corollary to Theorem \ref{thm:mainbij}]
        The number of $K$-Bianchi circles of curvature $\sqrt{-\Delta}f$, considered up equivalence, is $h_f$.  In particular, the number of $K$-Bianchi circles of curvature $\sqrt{-\Delta}f$ with centres in the fundamental parallelogram
        \[
                \{ a + b \omega : 0 \le a < 1, 0 \le b < 1 \}
        \]
        is $2h_f$, unless $f=1$ and $\Delta = -4$, in which case it is $1$.
\end{corollary}

\begin{proof}Follows immediately, except to observe that only in the case of $\Delta=-4$, $f=1$, do we have circles which reflect through the origin to themselves modulo $\OK$.
\end{proof}

We close the section with a brief remark on the relation between primeval lattices and the sequence \eqref{eqn:exact}.
Let $\mathcal{L}_f$ be the set of primeval lattices for the order $\Of$.  The units taking $\Lambda \in \mathcal{L}_f$ to itself are exactly $\Of^*$.  Therefore, the kernel of $\theta_f$ has a set of coset representatives given by $\mathcal{L}_f$ modulo $\mathcal{U}_f$.
On the other hand, the exact sequence
\eqref{eqn:exact}
gives another description of the kernel of $\theta_f$. 
In fact, we have a slightly stronger result:
\begin{proposition}
There is a bijection
\[
\mathcal{L}_f \rightarrow (\OK/f)^*/(\ZZ/f)^*
\]
The map is given by considering $\Lambda \in \mathcal{L}_f$ in $\OK$ and taking the indicated quotients.  The inverse is given by
\[
\beta \mapsto        {f}\OK + \beta \ZZ.
\]
\end{proposition}

\begin{proof}
        It follows from the fact that $\Lambda$ is a locally principal ideal that its image consists of just one\footnote{More directly, one could choose a Hermite basis for $\Lambda$ in terms of $1, \omega$; in this case one element of the basis is in $f\OK$, from which this observation follows.} element of $(\OK/f)^*/(\ZZ/f)^*$, so that the map is well-defined.  It is immediate to check that the maps are inverses (note that all primeval lattices have the form $f\OK + \beta\ZZ$).  
\end{proof}

Therefore, taking the quotient on either side by $\mathcal{U}_f$, we obtain a collection in bijection to $\SK/\sim$ via Theorem \ref{thm:mainbij}.

\section{Connectedness and Euclideanity}
\label{sec:euclidean}

The Schmidt arrangement is quite different in the cases of Euclidean and non-Euclidean rings of integers $\OK$.  Our purpose in this section is to prove the following.

\begin{theorem}
        \label{thm:connectedEuclidean}
        The Schmidt arrangement $\SK$ is connected if and only if $\OK$ is a Euclidean domain.  If $\SK$ is disconnected, it has infinitely many connected components.  
\end{theorem}

Write $E_2(\OK)$ for the subgroup of $\SL_2(\OK)$ generated by elementary matrices.  
That $\SK$ is connected for Euclidean domains is a consequence of a relationship between $E_2$ and Euclideanity due to Cohn.

\begin{theorem}[Cohn, {\cite{\Cohn}}]
        \label{thm:cohn}
 Let $K$ be a imaginary quadratic field.  Then $\OK$ is Euclidean if and only if $\SL_2(\OK)$ is generated by the elementary matrices.
\end{theorem}

For number fields besides imaginary quadratic fields, $\SL_2(\OK)$ is always generated by the elementary matrices \cite{MR0244257,MR0435293}.  Nica strengthened Cohn's result.

\begin{theorem}[Nica, {\cite{\Unreasonable}}]
        \label{thm:nica}
        Let $K$ be a imaginary quadratic field.  The subgroup $E_2(\OK)$ of $\SL_2(\OK)$ generated by the elementary matrices is either the whole group or else it is a non-normal infinite-index subgroup.
\end{theorem}

For more on this and surrounding results, the reader is also directed to the very nice exposition of \cite{\Unreasonable}.

\begin{definition}
        Define the \emph{tangency graph} of $\SK$ to be the graph whose vertices are the $K$-Bianchi circles of $\SK$ and whose edges represent tangencies.  Two circles are \emph{tangency-connected} if their vertices are in the same connected component of this graph.  We will refer to the subsets of $\SK$ corresponding to connected components of the tangency graph as \emph{tangency-connected components} of $\SK$.
\end{definition}

Cohn's theorem is enough to show that $\SK$ is tangency-connected if and only if $\OK$ is Euclidean, but tangency-connectedness is a stronger notion than connectedness.  To show that $\SK$ is disconnected when $\OK$ is not Euclidean, we do not use Cohn's result, and instead prove it directly.  As a consequence, this provides a new proof of one direction of Cohn and Nica's results:  \emph{If $\OK$ is not Euclidean, then the subgroup of $\PSL_2(\OK)$ generated by elementary matrices is of infinite index.}  That is, we will see below that the orbit of one $K$-Bianchi circle under elementary matrices lies in a single tangency-connected component of $\SK$, hence in one connected component of $\SK$.  However, we will show that there are infinitely many connected components.

In what follows, we will first discuss tangency-connectedness, which is fairly simple, before addressing the more difficult issue of connectedness.

\begin{proposition}
        \label{prop:tang-connect}
        The Schmidt arrangement $\SK$ is tangency-connected if and only if $\OK$ is a Euclidean domain.  If it is tangency-disconnected, it has infinitely many tangency-connected components.
\end{proposition}

\begin{proof}
        First, we claim that the following three matrices (multiplied on the right) take a circle to another in the same connected component of $\SK$:
        \[
                \begin{pmatrix} 1 & \tau \\ 0 & 1 \end{pmatrix}, 
                \begin{pmatrix} 1 & 1 \\ 0 & 1 \end{pmatrix}, 
                \begin{pmatrix} 0 & 1 \\ 1 & 0 \end{pmatrix}.
        \]
        The first is an example of moving to a tangent circle, according to Proposition \ref{prop:tangentfamilies}.  The second leaves the circle unchanged, and the third reverses orientation.
        But note that
        \[
                \begin{pmatrix} 1 & a \\ 0 & 1 \end{pmatrix}
                \begin{pmatrix} 1 & b \\ 0 & 1 \end{pmatrix} =
                \begin{pmatrix} 1 & a+b \\ 0 & 1 \end{pmatrix},
                \quad
                \begin{pmatrix} 1 & 0 \\ a & 1 \end{pmatrix} = 
                \begin{pmatrix} 0 & 1 \\ 1 & 0 \end{pmatrix}
                \begin{pmatrix} 1 & a \\ 0 & 1 \end{pmatrix}
                \begin{pmatrix} 0 & 1 \\ 1 & 0 \end{pmatrix}.
        \]
        Therefore, the orbit of $\widehat{\RR}$ under elementary matrices is contained in the tangency-connected component of $\widehat{\RR}$.  By Theorem \ref{thm:cohn}, then, we have shown that if $K$ is Euclidean, then $\SK$ is tangency-connected (and hence connected).

Now, suppose $K$ is not Euclidean (in particular, $K \neq \QQ(\sqrt{-3})$).

The circles tangent to $C$ are given by right multiplication by an elementary matrix (after changing orientation if necessary), so all circles which are in the tangency-connected component of $\widehat{\RR}$ are in the orbit of the image of $E_2(\OK)$ in $\PSL_2(\OK)$.  By Theorem \ref{thm:cohn}, $\SK$ is therefore tangency-disconnected.

The final statement is a consequence of Theorem \ref{thm:nica}.
\end{proof}

Tangency-disconnectedness is weaker than disconnectedness.  To prove the latter, we define a special circle which lies in the complement of the Schmidt arrangement.

\begin{figure}
        \includegraphics[width=4.5in]{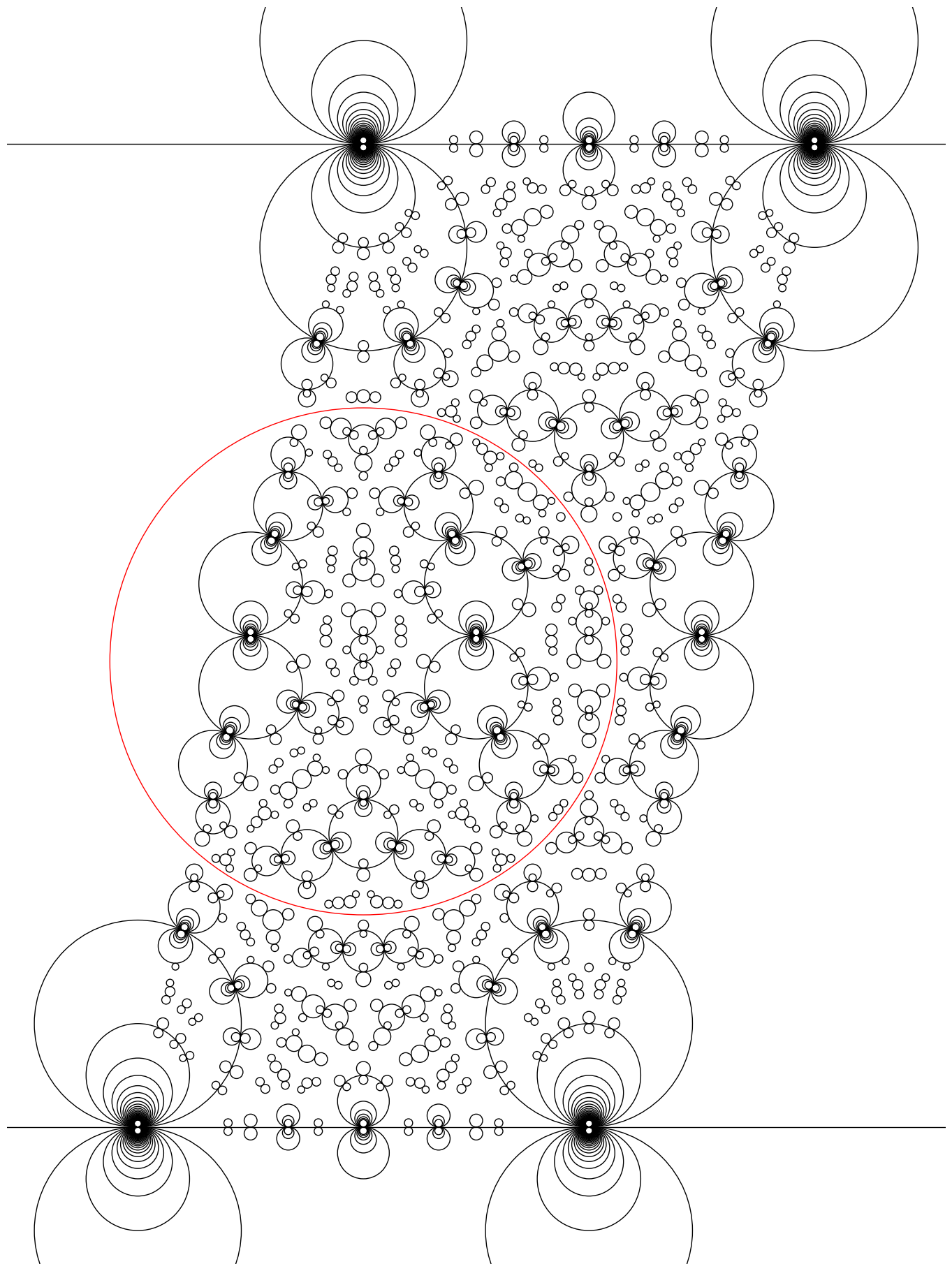}
        \caption{Schmidt arrangement of $\QQ(\sqrt{-19})$ (showing circles up to curvature $30\sqrt{19}$) with the ghost circle indicated in red.}
\label{fig:ghost}
\end{figure}

\begin{definition}
        The \emph{ghost circle} for a imaginary quadratic field $K$ is the circle orthogonal to the unit circle having center 
        \[
                \left\{ \begin{array}{ll}
                                \frac12 + \frac{\sqrt{\Delta}}{4} & \Delta \equiv 0 \pmod 4 \\
                         \frac12 + \frac{-\Delta-1}{4\sqrt{\Delta}} & \Delta \equiv 1 \pmod 4 \\
                        \end{array} \right. ,
                \]
                                and positive orientation.  (If no such circle exists, we say $K$ has no ghost circle.)
\end{definition}

\begin{lemma}
        If $\Delta > -11$, then $K$ has no ghost circle.  Otherwise, it has one, and its curvature is 
        \[
                \left\{ \begin{array}{ll}
                                4/\sqrt{-\Delta-12} & \Delta \equiv 0 \pmod 4 \\
      4\sqrt{-\Delta}/\sqrt{\Delta^2 +14\Delta + 1} & \Delta \equiv 1 \pmod 4 \\
                        \end{array} \right.
                \]
\end{lemma}

\begin{proof}
        Let $G$ be the ghost circle, of curvature $B$.  As it is orthogonal to the unit circle, its curvature and co-curvature are equal.  We use fact \eqref{eqn:cocurv}, which in the case $\Delta \equiv 0 \pmod 4$ has the form
        \[
                B^2 = \frac{B^2}4 + \frac{B^2\Delta}{16} - 1 \quad \implies \quad                B = \pm \frac{4}{\sqrt{-\Delta-12}}.
        \]
        and in the case $\Delta \equiv 1 \pmod 4$ has the form
\[
        B^2 = \frac{B^2}4 + \frac{B^2(-\Delta-1)^2}{16\Delta} - 1 \quad \implies \quad                B = \pm \frac{4\sqrt{-\Delta}}{\sqrt{\Delta^2 + 14\Delta+1}}.
        \]
        As $G$ is positively oriented, we choose the positive square root.  These roots are real for fundamental discriminants $< -11$ (the first such is $\Delta = -15$).
\end{proof}

\begin{lemma}
        \label{lemma:nogoghost}
        Suppose $K$ has a ghost circle.  Then no $K$-Bianchi circle intersects the ghost circle.
\end{lemma}

\begin{proof}
        We will show that the Pedoe product of the ghost circle $G$ with any $K$-Bianchi circle $C$ is greater than $1$ in absolute value.  By Proposition \ref{prop:pedoeproduct}, this implies that the circles do not intersect.  The proof is somewhat technical, relying on the fact that the image $\pi(C)$ is restricted in certain ways (for example, not all elements of $i\OK$ can occur as curvature-centres).
        
        Let $B$ be the curvature of $G$.  First, assume $\Delta \equiv 0 \pmod 4$.  Then 
\[
        \pi(G) =  \left(B, B, \frac{B}{2}, \frac{\sqrt{-\Delta} B}{4} \right)
        = \frac{B}{\sqrt{-\Delta}}\left( \sqrt{-\Delta}, \sqrt{-\Delta}, \frac{\sqrt{-\Delta}}{2}, \frac{-\Delta}{4} \right).
        \]
        By Proposition \ref{prop:curvature}, we know that
        \[
                \pi(C) \in \sqrt{-\Delta}\ZZ \times \sqrt{-\Delta}\ZZ
                \times \sqrt{-\Delta}\ZZ \times \ZZ.
        \]
        Hence, taking the product,
        \[
                \langle G, C \rangle \in \frac{B\sqrt{-\Delta}}{4}\ZZ.
                \]
                But 
                \[
                        \frac{B\sqrt{-\Delta}}{4} = \frac{\sqrt{-\Delta}}{\sqrt{-\Delta-12}} > 1.
                \]
                Hence it suffices to show that $\langle G, C \rangle \neq 0$.  Suppose $\pi(C) = (b,b',x,y)$.  By \eqref{eqn:cocurv},
                \[
                        bb' - x^2 - y^2 + 1 = 0
                \]
The first two terms are divisible by $\Delta$, hence $y$ is odd.  The equation $\langle G, C \rangle = 0$ is equivalent to
\[
        -4\frac{b+b'}{\sqrt{-\Delta}} + 2\frac{x}{\sqrt{-\Delta}} - y = 0
\]
which is impossible for odd $y$, as each of the displayed fractions is an integer.

The case $\Delta \equiv 1 \pmod 4$ is slightly more involved, but similar in spirit.  In this case,
\[
        \pi(G) = \frac{B}{\sqrt{-\Delta}} \left( \sqrt{-\Delta}, \sqrt{-\Delta}, \frac{\sqrt{-\Delta}}{2}, \frac{-\Delta-1}{4} \right).
\]
Suppose that $C$ has curvature-centre $z = i\frac{a+b\sqrt{\Delta}}{2}$, $a,b \in \ZZ$, curvature $c\sqrt{-\Delta}$ and co-curvature $c'\sqrt{-\Delta}$, where $z \in i\OK$ and $c,c' \in \ZZ$ by Proposition \ref{prop:curvature}.  By \eqref{eqn:cocurv}, $z\overline{z} = \Delta cc'+1 \equiv 1 \pmod \Delta$.  

With this notation, 
\[
        \pi(C) = \left( c\sqrt{-\Delta}, c'\sqrt{-\Delta}, -b\sqrt{-\Delta}/2, a/2 \right).
\]

Next, we show that $a \equiv 2 \pmod \Delta$.  By assumption, $z = i( \alpha \overline{\delta} - \gamma\overline{\beta} )$ (Proposition \ref{prop:curvature}), while $\alpha \delta - \gamma\beta = 1$.  Therefore $i - z = i(\alpha(\delta - \overline{\delta}) - \gamma( \beta - \overline{\beta}))$.  In particular,
\[
  i - \Im(z) \in \Delta \ZZ.
\]
This implies that $a \equiv 2 \pmod {2\Delta}$.

Let us write $a = a'\Delta + 2$, where $a' \equiv 0 \pmod 2$.  We must also have $b \equiv 0 \pmod 2$.
Now we compute
\begin{align*}
4\frac{\sqrt{-\Delta}}{B} \langle G, C \rangle 
&= 
4(c+c')\Delta + b\Delta - a \frac{\Delta + 1}{2} \\
&=
4(c+c')\Delta + b\Delta - \Delta a' \frac{\Delta + 1}{2} + (\Delta + 1). \notag
\end{align*}
This quantity lies in $2\Delta \ZZ + \Delta + 1$.  That is,
\[
        \frac{4}{B\sqrt{-\Delta}}\langle G, C \rangle = k + \frac{1}{\Delta},
        \]
        where $k \in 1 + 2\ZZ$.  To conclude that $|\langle G, C \rangle| > 1$, it suffices to show that 
        \begin{equation*}
                \frac{4}{B\sqrt{-\Delta}} < 1 - \frac{1}{\Delta}.
        \end{equation*}
        The inequality 
	is equivalent to
        \[
14\Delta^3 + 3\Delta^2 - 1 < 0
        \]
        which certainly holds in the range $\Delta < -11$, and this completes the proof.
\end{proof}

\begin{proof}[Proof of Theorem \ref{thm:connectedEuclidean}]
        Proposition \ref{prop:tang-connect} shows that if $\OK$ is Euclidean then $\SK$ is tangency-connected and therefore connected.  If $\OK$ is non-Euclidean, then $\Delta < -12$ and therefore there exists a ghost circle $G$ in the complement of $\SK$ by Lemma \ref{lemma:nogoghost}.  But since there is a circle of $\SK$ passing through every point $z \in K$, the circles of $\SK$ are dense and therefore populate both the interior and exterior of $G$.  However, they do not intersect $G$ itself.  But the interior and exterior are not connected, so that $\SK$ must be disconnected.

        Finally, any element of $\PSL_2(\OK)$ stabilizes $\SK$ and preserves intersections, therefore the orbit of $G$ under $\PSL_2(\OK)$ lies in the complement of $\SK$.  This implies that $\SK$ falls into infinitely many connected components (for example, take all translations of $G$ by $\OK$).
\end{proof}

{\bf Acknowledgements.} The author would like to thank Elena Fuchs, Robert Hines, Jerzy Kocik, Jeffrey Lagarias, Daniel Martin, Jonathan Wise and an anonymous referee 
for helpful discussions and suggestions.  The author would also like to thank Arseniy Sheydvasser for pointing out an error, since fixed, in the last theorem.  She is also indebted to the organisers, Alina Bucur, Pririta Paajanen and Lillian Pierce, of the 6th European Women in Mathematics Summer School on Apollonian Circle Packings.  Thank you to the lecturers Elena Fuchs and Hee Oh, to the EWM itself, and to the Institute Mittag-Leffler where the workshop was held in June of 2014.  The summer school provided the opportunity to learn a great deal.

\bibliography{app-pack-bib}{}
\bibliographystyle{plain}

\end{document}